\documentclass[11pt]{amsart}
\theoremstyle{plain}
\newtheorem{thm}{Theorem}[section]
\newtheorem{theorem}[thm]{Theorem}

\newtheorem{corollary}[thm]{Corollary}
\newtheorem{proposition}[thm]{Proposition}
\theoremstyle{definition}

\newtheorem{definition}[thm]{Definition}

\numberwithin{equation}{section}

\newcommand{\p}{\partial}

\newcommand{\C}{{\mathbb C}}

\newcommand{\BP}{{\mathbb P}}

\newcommand{\End}{{\rm End}}

\newcommand{\fg}{{\mathfrak g}}

\newcommand{\fgl}{{\mathfrak g}{\mathfrak l}}

\def\Sym{\mathop{\rm Sym}\nolimits}

\title[Symmetrizer group of a projective hypersurface]{Symmetrizer group of a projective hypersurface}

\author{Jun-Muk Hwang}

\thanks{This work was supported by the Institute for Basic Science (IBS-R032-D1).}
\begin{document}

\begin{abstract}
To each projective hypersurface which is not a cone, we associate an abelian linear algebraic group
called the symmetrizer group of the corresponding symmetric form. This group describes the set of homogeneous polynomials with the same Jacobian ideal and gives a conceptual explanation of  results by Ueda--Yoshinaga and Wang. In particular, the diagonalizable part of
the symmetrizer group detects Sebastiani-Thom property of the hypersurface and
its unipotent part  is related to the singularity of the hypersurface.

\end{abstract}

\maketitle



\section{Introduction}
We work over the complex numbers. Throughout, we fix an integer $d \geq 3$ and a vector space $V$ with $\dim V = n \geq 2$.  Let $\BP V$ be the projectivization of $V$, the set of 1-dimensional subspaces of $V$. For a nonzero vector $v \in V$, we denote by $[v] \in \BP V$ the 1-dimensional subspace $\C v \subset V$. We regard the vector space $\Sym^d V^*$ of symmetric $d$-forms on $V$ as the subspace of the vector space $\otimes^d V^*$ of $d$-linear forms on $V$, which  consists of all $ F \in \otimes^d V^*$ satisfying $$ F(v_1, \ldots, v_d) = F(v_{\sigma(1)}, \ldots, v_{\sigma(d)}) $$ for any permutation $\sigma$ of $\{1, \ldots, d\}$. Denote by $Z(F) \subset \BP V$ the hypersurface defined by the homogenous polynomial of degree $d$ corresponding to $F$.

\begin{definition} \begin{itemize} \item[(i)] For $F \in \Sym^d V^*$, define the homomorphism $\p F: V \to \Sym^{d-1} V^*$, called the {\em Jacobian} of $F$, by $$(\p F (u)) (v_1, \ldots, v_{d-1}) := F(u, v_1, \ldots, v_{d-1})$$ for all $u, v_1, \ldots, v_{d-1} \in V.$  \item[(ii)] We say that $F$ is {\em nondegenerate} if ${\rm Ker}(\p F) =0,$ equivalently, if the hypersurface $Z(F)$ is not a cone. Denote by $\Sym^d_o V^* $ the Zariski-open subset of $\Sym^d V^*$ consisting of nondegenerate forms.  \item[(iii)] Let ${\rm Gr}(n; \Sym^{d-1}V^*)$ be the Grassmannian of $n$-dimensional subspaces in the vector space $\Sym^{d-1} V^*.$
For a nondegenerate form $F \in \Sym^d_o V^*$, the dimension of ${\rm Im}(\p F)$ is $n$. Let $J(F) $ be the point in ${\rm Gr}(n; \Sym^{d-1}V^*)$ corresponding to the $n$-dimensional subspace ${\rm Im}(\p F)$. This defines a morphism $$J: \Sym^d_o V^* \to {\rm Gr}(n; \Sym^{d-1}V^*).$$ \end{itemize} \end{definition}

Historically, the interest in the  morphism $J$ arose from the study of variation of Hodge structures for families of hypersurfaces. In particular, Carlson and Griffiths showed in Section 4.b of  \cite{CG} that $J^{-1}(J(F)) = \C^{\times}\cdot F$ for a general $F\in \Sym^d_o V^*$. Obvious examples satisfying $J^{-1}(J(F)) \neq \C^{\times} \cdot F$ are symmetric forms of Sebastiani-Thom type in the following sense.

\begin{definition}\label{d.ST}
An element $F \in \Sym^d V^*$ is {\em of Sebastiani-Thom type}, if there is a direct sum decomposition $V = V_1 \oplus \cdots \oplus V_k$ with $ k \geq 2$ such that $F= F_1 + \cdots + F_k$ for some $F_i \in \Sym^d V^*_i, 1 \leq i \leq k.$ More precisely, this means that for each $1 \leq i \neq j \leq k$ and  $u \in V_i, w \in V_j,$  $$F_i(u, w, v_1, \ldots, v_{d-2}) =0$$ for all $v_1, \ldots, v_{d-2} \in V.$
\end{definition}

For $F = F_1 + \cdots + F_k \in \Sym^d_o V^*$ of Sebastiani-Thom type, we have $$J(F) = J(c_1 F_1 + \cdots + c_k F_k)$$ for any $c_1, \ldots, c_k \in \C^{\times}$. Thus $J^{-1}(J(F)) \neq \C^{\times} \cdot F$.  Ueda and Yoshinaga proved the following in \cite{UY}.

\begin{theorem}\label{t.UY} Let $F \in \Sym^d_o V^*$ be such that   the  hypersurface $Z(F) \subset \BP V$ is nonsingular. Then $J^{-1}(J(F)) \neq \C^{\times} \cdot F$ if and only if $F$ is of Sebastiani-Thom type. \end{theorem}

By this result,  the key issue in understanding symmetric forms satisfying  $J^{-1}(J(F)) \neq \C^{\times} \cdot F$ is to study its relation with the singularity of the hypersurface $Z(F)$.  In this direction, Wang proved in \cite{Wa} the following result.

\begin{theorem}\label{t.Wang}
Let $F \in \Sym^d_o V^*$ be such that  $J^{-1}(J(F)) \neq \C^{\times} \cdot F$ and $F$ is not of Sebastiani-Thom type. Then the hypersurface $Z(F)$ has a singular point of multiplicity $d-1$, namely, there is a nonzero vector $u \in V$ such that $$F(u, u, v_1, \ldots, v_{d-2}) =0 \mbox{ for all } v_1, \ldots v_{d-2} \in V.$$ \end{theorem}

Of course, Theorem \ref{t.Wang} is reduced to Theorem \ref{t.UY} when $d=3$, but for $d \geq 4$, it gives additional information on the singularity of $Z(F)$.

The proofs of Theorem \ref{t.UY} and Theorem \ref{t.Wang} in \cite{UY} and \cite{Wa} are computational. Our main result is  a geometric description of the fibers of the morphism $J$, which gives  a more conceptual explanation of Theorems \ref{t.UY} and \ref{t.Wang}. More precisely, we describe the fibers of $J$ as follows.

\begin{theorem}\label{t.main}
Let $x \in {\rm Im}(J) \subset {\rm Gr}(n; \Sym^{d-1}V^*)$ be a point in the image of the morphism $J.$  Then there is a connected abelian algebraic subgroup $G_x \subset {\rm GL}(V)$ canonically associated to $x$, which contains $\C^{\times} \cdot {\rm Id}_V$,
such that the fiber $J^{-1}(x)$ is a principal homogeneous space of the group $G_x$.
\end{theorem}

By the classification of connected abelian groups (Theorems 3.1.1 and 3.4.7 of \cite{Sp}), we have a decomposition into direct product of algebraic groups $$G_x/(\C^{\times} \cdot {\rm Id}_V) \ = \ G_x^{\times} \times G_x^+, $$ where $G_x^{\times}$ is a diagonalizable group (an algebraic torus)  and $G_x^+$ is a vector group.
We have the corresponding decomposition of the Lie algebra $$ \fg_x / (\C \cdot {\rm Id}_V) \ = \  \fg_x^{\times} \oplus \fg_x^+.$$
We prove the following,  which is a refinement of Theorems \ref{t.UY} and \ref{t.Wang}.

\begin{theorem}\label{t.alpha} In the setting of Theorem \ref{t.main},
\begin{itemize} \item[(i)] any element  $F \in J^{-1}(x)$ is of Sebastiani-Thom type if and only if $\fg_x^{\times} \neq 0$; and \item[(ii)] if $\fg_x^+ \neq 0$, then  there is $0 \neq u \in V$ such that $[u]$ is a point of multiplicity $d-1$ on the hypersurface $Z(F)$ for all $F \in J^{-1}(x)$.  \end{itemize} \end{theorem}

 The diagonalizable part $G_x^{\times}$ of the group $G_x$ is well-explained by Theorem \ref{t.alpha} (i), but the unipotent part $G_x^+$  is not fully described by  (ii).
 It would be interesting to find more geometric consequences of the unipotent part $G^+_x$. We obtain the following result in this direction.

\begin{theorem}\label{t.count} In the setting of Theorem \ref{t.main}, assume that for some $F \in J^{-1}(x)$, the hypersurface $Z(F) \subset \BP V$ has only finitely many singular points of multiplicity $d-2$.
(This is the case, for example, if $Z(F)$ has only isolated singularities.)
\begin{itemize} \item[(i)] The number of points in $\BP \fg_x^+$ corresponding to nonzero elements $h \in \End(V)$ satisfying $h^2=0$ is less than or equal to the number of singular points of multiplicity $d-2$ on $Z(F)$.
\item[(ii)] For any $f \in \fg^+_x$, we have $f^3 =0.$ \end{itemize} \end{theorem}

 It would be interesting to investigate how big $\dim \fg_x^+$ can be, especially in the setting of Theorem \ref{t.count}.

The proof of Theorem \ref{t.main} is rather simple, once one realizes what the group $G_x$ should be.
We describe the group $G_x$  and prove Theorem \ref{t.main}  in  Section \ref{s.group}. Theorems \ref{t.alpha} and \ref{t.count} are proved in Section \ref{s.proof}.

\section{Symmetrizer group of a symmetric form}\label{s.group}
\begin{definition}\label{d.algebra}
For $F \in \Sym^d V^*$ and $g \in \End(V)$, define $F^g \in \otimes^d V^*$ by
$$F^g(v_1, \ldots, v_d) := F(g \cdot v_1, v_2, \ldots, v_d) \mbox{ for } v_1, \ldots, v_d \in V.$$ We say that $g$ is a {\em symmetrizer} of $F$ if $F^g \in \Sym^d V^*$, namely, \begin{eqnarray*} F(g \cdot v_1, v_2, v_3, \ldots, v_d) & = & F(v_1, g \cdot v_2, v_3, \ldots, v_d)\\ & = & F(v_1, v_2, g \cdot v_3, \ldots, v_d) \\ &= & \cdots \\ & = & F(v_1, v_2, v_3, \ldots, g \cdot v_d). \end{eqnarray*} The subspace $\fg_F \subset \End(V)$  of all symmetrizers of $F$ is called the {\em symmetrizer algebra} of $F$ and the intersection $G_F := \fg_F \cap {\rm GL}(V)$ is called the {\em symmetrizer group} of $F$. \end{definition}

The two names, the symmetrizer algebra and the symmetrizer group, in Definition \ref{d.algebra} are justified by the next two propositions.

\begin{proposition}\label{p.abelian} In Definition \ref{d.algebra}, the following holds.
\begin{itemize} \item[(i)] If $g, h \in \fg_F$, then $g \circ h \in \fg_F$. In particular, the vector space $\fg_F$ is a subalgebra under the composition in $\End(V)$, hence a  Lie subalgebra of $\End(V).$
\item[(ii)] If $F$ is nondegenerate, then $g \circ h = h \circ g$ for any $g, h \in \fg_F$, namely, the Lie algebra $\fg_F$ is abelian.
    \item[(iii)] If $h \in \fg_F$, then $F(u, w, v_1, \ldots, v_{d-2}) =0$ for any $u \in {\rm Im}(h), w \in {\rm Ker}(h)$ and $v_1, \ldots, v_{d-2} \in V.$ \end{itemize} \end{proposition}

    \begin{proof}
    Write $gh = g \circ h$ for simplicity. Then for $g,h \in \fg_F$, \begin{eqnarray*}
    F(gh \cdot v_1, v_2, v_3, \ldots, v_d) & = & F(h \cdot v_1, v_2, g \cdot v_3, \ldots, v_d) \\ &= & F(v_1, h \cdot v_2, g \cdot v_3, \ldots, v_d) \\ & =&
    F(v_1, gh \cdot v_2, v_3, \ldots, v_d). \end{eqnarray*}
    Thus $gh \in \fg_F$,  proving (i).

    Now assume that $F$ is nondegenerate. Then \begin{eqnarray*}
    F(gh \cdot v_1, v_2, v_3, \ldots, v_d) & = & F(h \cdot v_1, g \cdot v_2,  v_3, \ldots, v_d) \\ &= & F(v_1, hg \cdot v_2,  v_3, \ldots, v_d) \\ &=& F(hg \cdot v_1, v_2, v_3, \ldots, v_d),  \end{eqnarray*} where the last equality uses (i).
Thus $F((hg - gh) \cdot v_1, v_2, \ldots, v_d) =0$ for all $v_1, \ldots, v_d \in V$.
By the nondegeneracy of $F$, this implies $hg= gh$, proving (ii).

To prove (iii), write $u = h\cdot u'$ for some $u' \in V$. Then for $w \in {\rm Ker}(h)$,
\begin{eqnarray*} F(u, w, v_1, \ldots, v_{d-2}) &=& F(h \cdot u', w, v_1, \ldots, v_{d-2}) \\ &=& F(u', h \cdot w, v_1, \ldots, v_{d-2}) = 0.\end{eqnarray*}
\end{proof}

\begin{proposition}\label{p.group}
For $F \in \Sym^d V^*$, the intersection $G_F=\fg_F \cap {\rm GL}(V)$
is a connected subgroup of ${\rm GL}(V),$ corresponding to the Lie subalgebra $\fg_F \subset {\rm End}(V) = \fgl(V).$  It is an abelian group if $F$ is nondegenerate. \end{proposition}

\begin{proof}
Since $G_F$ is a Zariski open subset in the vector space $\fg_F,$ it is connected.
To check that it is a subgroup, it suffices to show,  by Proposition \ref{p.abelian}, that  $g \in G_F$ implies $g^{-1} \in G_F$.  For $v_i \in V$, write $u_i := g^{-1} \cdot v_i$. Then
\begin{eqnarray*} F(g^{-1} \cdot v_1, v_2, v_3, \ldots, v_d) & = & F(u_1, g \cdot u_2, v_3, \ldots, v_d) \\ & = & F(g \cdot u_1, u_2, v_3, \ldots, v_d) \\ & = & F(v_1, g^{-1} \cdot v_2, v_3, \ldots, v_d). \end{eqnarray*} This shows that $g^{-1} \in G_F$. \end{proof}

The proof of the following proposition is straightforward.

\begin{proposition}\label{p.ST}
When $F = F_1 + \cdots + F_k \in \Sym^d V^*$ is of Sebastiani-Thom type in the notation of Definition \ref{d.ST}, $$\fg_F = \fg_{F_1} \oplus \cdots \oplus \fg_{F_k} \mbox{ and } G_F = G_{F_1} \times \cdots \times G_{F_k},$$ where the products mean those coming from the decomposition $V= V_1 \oplus \cdots \oplus V_k$. \end{proposition}

\begin{proposition}\label{p.basic}
For $F \in \Sym^d V^*$ and $g \in G_F$,
\begin{itemize} \item[(i)] $G_{F^g} = G_F;$
\item[(ii)] ${\rm Ker}(\p F^g) = g^{-1} \cdot {\rm Ker}(\p F);$ and
\item[(iii)] $F^g$ is nondegenerate if and only if $F$ is nondegenerate.
\end{itemize}
Assume that $F$ is nondegenerate, then
\begin{itemize}
\item[(iv)] $F = F^g$ if and only if $g= {\rm Id}_V$; and
\item[(v)] $J(F) = J(F^g)$. \end{itemize} \end{proposition}

\begin{proof} To prove (i), pick  $h \in G_F$. Then $$F^g (h \cdot v_1, \ldots, v_d) = F(gh \cdot v_1, \ldots, v_d)$$ is symmetric in $v_1, \ldots, v_d$ because $gh \in G_F$. Thus $h \in G_{F^g}$, proving $G_F \subset G_{F^g}$. In particular,  if $f \in G_{F^g}$, then $g^{-1}f \in G_{F^g}$. Consequently,  $$F( f \cdot v_1, \ldots, v_d) = F^g (g^{-1}f \cdot v_1, \ldots, v_d)$$ is symmetric in $v_1, \ldots, v_d$. This shows $f \in G_F$, proving $G_{F^g} \subset G_F$.

Note that $u \in {\rm Ker}(\p F^g)$ if and only if $$F^g(u, v_1, \ldots, v_{d-1})= F(g \cdot u, v_1, \ldots, v_{d-1}) =0 \mbox{ for all } v_1, \ldots, v_{d-1} \in V.$$ This is equivalent to saying $g \cdot u \in {\rm Ker}(\p F)$. This proves (ii). (iii) is immediate from (ii).

Now assume that $F$ is nondegenerate and $F = F^g$ for some $g \in G_F$.
Then $$0= F(v_1, v_2, \ldots, v_{d}) - F(g \cdot v_1, v_2, \ldots, v_{d}) = F(({\rm Id}_V - g ) \cdot v_1, v_2, \ldots, v_d)$$ for all $v_1,  v_2, \ldots, v_{d-1} \in V$.
By the nondegeneracy of $F$, this implies $g = {\rm Id}_V$, proving (iv).

To check (v), for each $u \in V$ and $v_1, \ldots, v_{d-1} \in V$, $$(\p F(u)) (v_1, \ldots, v_{d-1}) = F(u, v_1, \ldots, v_{d-1}) = F^g(g^{-1}\cdot u, v_1, \ldots, v_{d-1}).$$ This means $\p F(u) = \p F^g(g^{-1} \cdot u).$ Thus ${\rm Im}(\p F) = {\rm Im}(\p F^g)$, implying $J(F) = J(F^g).$ \end{proof}

The following is the converse of Proposition \ref{p.basic} (iv).

\begin{proposition}\label{p.converse}
Let $F, \widetilde{F} \in \Sym^d V^*$ be nondegenerate  symmetric forms satisfying $J(F) = J(\widetilde{F})$. Then there exists $g \in G_F$ such that $\widetilde{F} = F^g.$ \end{proposition}

\begin{proof} Let $g \in {\rm GL}(V) $ be the composite $$V \stackrel{\p \widetilde{F}}{\longrightarrow} {\rm Im}(\p \widetilde{F}) = {\rm Im}(\p F) \stackrel{(\p F)^{-1}}{\longrightarrow} V.$$ Then $g \in G_F$ and $F^g = \widetilde{F}$ because \begin{eqnarray*} F^g(v_1, v_2, \ldots, v_d) &=& F( g \cdot v_1, v_2, \ldots, v_d) \\ & =& (\p F (g \cdot v_1)) (v_2, \ldots, v_d) \\ & =& (\p F \circ g (v_1)) ( v_2, \ldots, v_d) \\ &=& (\p \widetilde{F} (v_1)) (v_2, \ldots, v_d) \\ &=& \widetilde{F}( v_1, v_2, \ldots, v_d) \end{eqnarray*} for all $v_1, \ldots, v_d \in V.$
 \end{proof}

The following direct corollary of Propositions \ref{p.basic} and \ref{p.converse} implies  Theorem \ref{t.main}.

\begin{corollary}\label{c.Gx}
For each point $x \in {\rm Im}(J) \subset {\rm Gr}(n; \Sym^{d-1}V^*),$
define $G_x := G_F \subset {\rm GL}(V)$ for any  $F \in  J^{-1}(x)$. Then
 $G_x$ does not depend on the choice of $F \in J^{-1}(x)$ and the fiber $J^{-1}(x)$ is a principal homogeneous space of $G_x$.
  \end{corollary}

\section{Diagonalizable and unipotent components of the symmetrizer group}\label{s.proof}

In this section, we fix a nondegenerate form $F \in \Sym^d V^*$. The  connected abelian group $G_F/ (\C^{\times} \cdot {\rm Id}_V)$ has  a canonical decomposition $$  G_F/(\C^{\times} \cdot {\rm Id}_V)= G_F^{\times} \times G_F^+,$$ where $G_F^{\times}$ is an algebraic torus and $G_F^+$ is a vector group. Let $$  \fg_F/(\C \cdot {\rm Id}_V) = \fg_F^{\times} \oplus \fg_F^+$$ be the corresponding decomposition of the Lie algebra, where $\fg_F^{\times}$ (resp. $\fg_F^+$) consists of semi-simple (resp.  nilpotent) elements.
The next proposition implies Theorem \ref{t.alpha} (i).

\begin{proposition}\label{p.mult}
If $\fg_F^{\times} \neq 0$,  then $F$ is of Sebastiani-Thom type. More precisely, let $$V = V_1 \oplus \cdots \oplus V_k, \ k \geq 2$$ be the decomposition into distinct weight spaces of the diagonalizable subgroup $\widetilde{G}_F^{\times} \subset G_F \subset {\rm GL}(V),$ which is the inverse image of the diagonalizable subgroup $G_F^{\times}$. Then $F = F_1 + \cdots + F_k$ for some $F_i \in \Sym^d V^*_i$. \end{proposition}

\begin{proof}
 We claim that if $v_i \in V_i$ and $v_j \in V_j$ for $i \neq j$, then  $$F(v_i, v_j, u_1, \ldots, u_{d-2}) =0 \mbox{ for any } u_1, \ldots, u_{d-2} \in V.$$ By the claim, if we set $F_i = F|_{V_i}$, then we obtain $F= F_1 + \cdots + F_k.$

To prove the claim, let $\lambda_i$ be the weight of $V_i$ for $1\leq i \leq k$.
Then for any $g \in \widetilde{G}_F^{\times}$,
\begin{eqnarray*} \lambda_i(g) F(v_i, v_j, u_1, \ldots, u_d) & = & F(g \cdot v_i, v_j, u_1, \ldots, u_{d-2}) \\ &=& F(v_i, g \cdot v_j, u_1, \ldots, u_{d-2}) \\ & =& \lambda_j(g) F(v_i, v_j, u_1, \ldots, u_{d-2}). \end{eqnarray*} Since $\lambda_i \neq \lambda_j$, the claim follows. \end{proof}

Propositions \ref{p.ST} and \ref{p.mult} show that the genuinely interesting part of the symmetrizer group is $G^+_F$.
The next proposition implies Theorem \ref{t.alpha} (ii).

\begin{proposition}\label{p.square}
Assume that $\fg^+_F \neq 0$. \begin{itemize}
\item[(i)] There exists at least one nonzero element $h \in \fg_F^+$ satisfying $h^2 =0$.
\item[(ii)] For $0 \neq h \in \fg^+_F$ satisfying $h^2=0$,  every  point of  $\BP {\rm Im}(h) \subset \BP V$ is a singular point of $Z(F)$ with  multiplicity $d-2$.  \end{itemize} \end{proposition}

\begin{proof}
Pick $0 \neq f \in \fg^+_F.$ All  powers $f^k, k \geq 2,$  belong to $\fg^+_F$ by Proposition \ref{p.abelian}. Since $f$ is nilpotent, there is an integer $\ell >1$ satisfying $f^{\ell} =0.$ Then $h := f^{\ell -1}$ satisfies $h^2 =0$. This proves (i).

In (ii),
from Proposition \ref{p.abelian} (iii) and  ${\rm Im}(h) \subset {\rm Ker}(h)$, we have $$F(u, u, v_1, \ldots, v_{d-2}) =0$$ for all $u \in {\rm Im}(h)$ and $v_1, \ldots, v_{d-2} \in V$. Thus $[u] \in \BP V$ is a singular point of $Z(F)$ with multiplicity $d-2$. \end{proof}

We reformulate Theorem \ref{t.count} as follows.

\begin{theorem}\label{t.count2}
Assume that $Z(F)$ has only finitely many singular points of multiplicity $d-2$.
\begin{itemize} \item[(i)]  If $0 \neq h \in \fg^+_F$ satisfies $h^2 =0,$ then $\dim {\rm Im}(h) =1.$
\item[(ii)] If $h, \widetilde{h} \in \fg^+_F$ are nonzero elements satisfying $h^2 = \widetilde{h}^2 =0$ and ${\rm Im}(h) = {\rm Im}(\widetilde{h})$, then $[h] = [\widetilde{h}] \in \BP \fg^+_F.$
\item[(iii)]
The number of points in $\BP \fg^+_F$ corresponding to nonzero elements $h \in \fg^+_F$ satisfying $h^2=0$ is less than or equal to the number of singular points of multiplicity $d-2$ on $Z(F)$.
\item[(iv)] For any $f \in \fg^+_F$, we have $f^3 =0.$ \end{itemize} \end{theorem}

    \begin{proof}
    By Proposition \ref{p.square}, if $h^2 =0$, then every point on $\BP {\rm Im}(h)$ is a singular point of $Z(F)$ with multiplicity $d-2$. By the assumption that there are only finitely many singular points of multiplicity $d-2$, we see $\dim {\rm Im}(h) =1$, proving (i).

To prove (ii), pick a nonzero element $u \in {\rm Im}(h) = {\rm Im}(\widetilde{h})$. By Proposition \ref{p.abelian} (iii), for all $ v_1, \ldots, v_{d-2} \in V,$
    $$F(u, {\rm Ker}(h), v_1, \ldots, v_{d-2}) =0 = F(u, {\rm Ker}(\widetilde{h}), v_1, \ldots, v_{d-2}).$$ If  $V = {\rm Ker}(h) + {\rm Ker}(\widetilde{h}),$ we have  $F(u, V, v_1, \ldots, v_{d-2}) =0,$
 a contradiction to the nondegeneracy of $F$. Thus $V \neq {\rm Ker}(h) + {\rm Ker}(\widetilde{h}),$ which implies ${\rm Ker}(h) = {\rm Ker}(\widetilde{h})$. It follows that $[h] = [\widetilde{h}]$.

(iii) follows from (ii), because $\BP {\rm Im}(h) \in \BP V$ is a singular point of $Z(F)$ with multiplicity $d-2$ by Proposition \ref{p.square} (ii).

    To prove (iv), assume the contrary  that for some $f \in \fg^+_F$ and an integer  $\ell \geq 4,$ the elements  $f, f^2, \ldots, f^{\ell-1}$  are nonzero and $f^{\ell} =0.$
    Then $h:= f^{\ell -1}$ and $\widetilde{h}:= f^{\ell -2}$ satisfy $h^2 = \widetilde{h}^2 =0$. Since ${\rm Im}(h) \subset {\rm Im}(\widetilde{h})$, we see by (i) and (ii) that $\widetilde{h} = c h$ for some $c \in \C^{\times}$. In other words, we have $f^{\ell-2} = c f^{\ell -1}$. Hence $f^{\ell-1} = c f^{\ell} =0$, a contradiction. \end{proof}

\bigskip
{\bf Acknowledgment} I would like to thank Baohua Fu, Keiji Oguiso and Kazushi Ueda for encouragement.

\bigskip
Institute for Basic Science

Center for Complex Geometry

Daejeon, 34126

Republic of Korea

jmhwang@ibs.re.kr

\begin{thebibliography}{6}
\bibitem{CG} Carlson, J., Griffiths, P.: Infinitesimal variations of Hodge structure and the global Torelli
problem. In: {\em Journées de Géometrie Algébrique d’Angers},
pp. 51–76. Sijthoff and Noordhoff (1980)

\bibitem{Sp} Springer, T. A.: {\em Linear algebraic groups}, Second Edition,
Birkh\"auser (1998)

\bibitem{UY}
Ueda, K.; Yoshinaga, M.: Logarithmic vector fields along smooth divisors in projective spaces. Hokkaido Math. J. 38 (2009)  409–415

\bibitem{Wa}
Wang, Z.: On homogeneous polynomials determined by their Jacobian ideal. Manuscripta Math. 146 (2015)  559–574
     \end{thebibliography}
\end{document}